\documentclass{amsart}

\usepackage{amsmath,amssymb,amsfonts}
\usepackage{stmaryrd}
\usepackage[margin=1in]{geometry}
\usepackage{url}
\usepackage{xcolor}
\usepackage{listings}

\lstdefinelanguage{GAP}{
  morekeywords={F, rels, PSL2Z, LowIndexSubgroupsFpGroup, Length, Print, for, if, then, else, fi, do, od},
  sensitive=true,
  morecomment=[l]{\#},
  morestring=[b]"
}
\newtheorem{theorem}{Theorem}[section]
\newtheorem{lemma}[theorem]{Lemma}
\newtheorem{cor}[theorem]{Corollary}
\newtheorem{prop}[theorem]{Proposition}

\theoremstyle{definition}
\newtheorem{definition}[theorem]{Definition}

\theoremstyle{remark}

\numberwithin{equation}{section}

\newcommand{\FF}{\mathbb{F}}

\newcommand{\gN}{\mathcal{N}}

\newcommand{\bC}{\mathbf{C}}
\newcommand{\bZ}{\mathbf{Z}}
\newcommand{\bF}{\mathbf{F}}
\newcommand{\bN}{\mathbf{N}}

\newcommand{\Aut}{\operatorname{Aut}}

\newcommand{\Irr}{\operatorname{Irr}}

\newcommand{\Syl}{\operatorname{Syl}}

\newcommand{\Sp}{\operatorname{Sp}}

\newcommand{\Gal}{\operatorname{Gal}}

\newcommand{\CHA}{\operatorname{char}}

\newcommand{\GF}{\mathrm{GF}}

\newcommand{\stab}{\operatorname{stab}}

\begin{document}

\title{Large orbits of nilpotent subgroups of linear groups}

\author{Yuchen Xu and Yong Yang}

\address{The Lawrenceville School, Lawrenceville, NJ 08648, USA}
\email{yuchennnxu@gmail.com}
\address{Department of Mathematics, Texas State University, 601 University Drive, San Marcos, TX 78666, USA.}


\makeatletter
\email{yang@txstate.edu}
\makeatother

\Large

\date{}



\begin{abstract}
Suppose that $G$ is a finite solvable group and $V$ is a finite, faithful and completely reducible $G$-module. Let $N$ be a nilpotent subgroup of $G$, then there exits $v \in V$ such that $|\bC_N(v)| \leq (|N|/p)^{1/p}$, where $p$ is the smallest prime divisor of $|N|$. 
\end{abstract}

\maketitle
\section{Introduction} \label{sec:introduction8}
   Let $G$ be a finite group and $V$ a finite, faithful and completely reducible $G$-module. It is a classical theme to study orbit structures of $G$ acting on $V$. One of the most important and natural questions about orbit structure is to establish the existence of an orbit of a certain size. For a long time, there has been a deep interest and need to examine the size of the largest possible orbits in linear group actions. The orbit $\{v^g \ |\ g \in G\}$ is called regular, if $\bC_G(v)=1$ holds or equivalently the size of the orbit $v^G$ is $|G|$.


Isaacs proved the following result ~\cite[Theorem A]{IMI2}. Let $N$ be a nontrivial $p$-group that acts faithfully on a group $H$, where $|H|$ is not divisible by $p$. Then there exists an element $x \in H$ such that $|\bC_N(x)| \leq {(|N|/p)}^{1/p}$ where $p$ is the smallest prime divisor of $|N|$. By the Hartley-Turull's Lemma ~\cite[Lemma 2.6.2]{HTurull}, this could be reduced to the following statement. Let $P$ be a nontrivial $p$-group that acts faithfully on a vector space $V$, where $|V|$ is not divisible by $p$. Then there exists an element $x \in V$ such that $|\bC_P(x)| \leq (|P|/p)^{1/p}$. Our result generalizes the case of the linear group actions and show the same bound hold for an arbitrary $p$-subgroup of a solvable linear group where the action is completely reducible. This also can be used to strengthen a result of Keller and the author ~\cite[Theorem 1.2]{KellerYang}.

The technique of the proof was developed in some recent papers of the second author ~\cite{Yang19, YYGluck}. The method, in some cases, can provide a good estimation on the size of the orbit. We remark that Isaacs's proof was extremely elegant and he found a way to avoid the detailed analysis of primitive linear groups. Though the proof in this paper is mainly calculation based, we were able to get something slightly stronger.

\section{Notation and Lemmas} \label{sec:Notation and Lemmas}

Notation:
\begin{enumerate}




\item If $V$ is a finite vector space of dimension $n$ over $\GF(q)$, where $q$ is a prime power, we denote by $\Gamma(q^n)=\Gamma(V)$ the semilinear group of $V$, i.e.,
\[\Gamma(q^n)=\{x \mapsto ax^{\sigma}\ |\ x \in \GF(q^n), a \in \GF(q^n)^{\times}, \sigma \in \Gal(\GF(q^n)/\GF(q))\},\] and we define \[\Gamma_0(q^n)=\{x \mapsto ax\ | \ x \in \GF(q^n), a \in \GF(q^n)^{\times}\}.\]\item We use $H \wr S$ to denote the wreath product of $H$ with $S$ where $H$ is a group and $S$ is a permutation group.
\end{enumerate}

\begin{definition} \label{defineEi}
Suppose that a finite solvable group $G$ acts faithfully, irreducibly and quasi-primitively on a finite vector space $V$. Let $\bF(G)$ be the Fitting subgroup of $G$ and $\bF(G)=\prod_i P_i$, $i=1, \dots, m$ where $P_i$ are normal $p_i$-subgroups of $G$ for different primes $p_i$. Let $Z_i = \Omega_1(\bZ(P_i))$. We define \[E_i=\left\{ \begin{array}{lll} \Omega_1(P_i) & \mbox{if $p_i$ is odd}; \\ \lbrack P_i,G,\cdots, G \rbrack & \mbox{if $p_i=2$ and $\lbrack P_i,G,\cdots, G \rbrack \neq 1$}; \\  Z_i & \mbox{otherwise}. \end{array} \right.\] By proper reordering we may assume that $E_i \neq Z_i$ for $i=1, \dots, s$, $0 \leq s \leq m$ and $E_i=Z_i$ for $i=s+1, \dots, m$. We define $E=\prod_{i=1}^s E_i$, $Z=\prod_{i=1}^s Z_i$ and we define $\bar{E}_i=E_i/Z_i$, $\bar{E}=E/Z$. Furthermore, we define $e_i=\sqrt {|E_i/Z_i|}$ for $i=1, \dots, s$ and $e=\sqrt{|E/Z|}$.
\end{definition}

\begin{theorem} \label{Strofprimitive}

Suppose that a finite solvable group $G$ acts faithfully, irreducibly and quasi-primitively on an $n$-dimensional finite vector space $V$ over finite field $\FF$ of characteristic $r$. We use the notation in Definition ~\ref{defineEi}. Then every normal abelian subgroup of $G$ is cyclic and $G$ has normal subgroups $Z \leq U \leq F \leq A \leq G$ such that,
\begin{enumerate}
\item $F=EU$ is a central product where $Z=E \cap U=\bZ(E)$ and $\bC_G(F) \leq F$;
\item $F/U \cong E/Z$ is a direct sum of completely reducible $G/F$-modules;
\item $E_i$ is an extra-special $p_i$-group for $i=1,\dots,s$ and $e_i=p_i^{n_i}$ for some $n_i \geq 1$. Furthermore $(e_i,e_j)=1$ when $i \neq j$ and $e=e_1 \dots e_s$ divides $n$, also $\gcd(r,e)=1$;
\item $A=\bC_G(U)$ and $G/A \lesssim \Aut(U)$, $A/F$ acts faithfully on $E/Z$;
\item $A/\bC_A(E_i/Z_i) \lesssim \Sp(2n_i,p_i)$;
\item $U$ is cyclic and acts fixed point freely on $W$ where $W$ is an irreducible submodule of $V_U$;
\item $|V|=|W|^{eb}$ for some integer $b$ and $|G:A| \mid \dim(W)$.
\item $G/A$ is cyclic.
\end{enumerate}
\end{theorem}
\begin{proof}
This is ~\cite[Theorem 2.2]{YY2}.
\end{proof}

\begin{theorem} \label{quote1}
Suppose that a finite solvable group $G$ acts faithfully, irreducibly and quasi-primitively on a finite vector space $V$. By Definition ~\ref{defineEi} and Theorem ~\ref{Strofprimitive}, $G$ will have a uniquely determined normal subgroup $E$ which is a direct product of extra-special $p$-groups for various $p$ and $e=\sqrt{|E/\bZ(E)|}$. Assume $e=5,6,7$ or $e \geq 10$ and $e \neq 16$, then $G$ will have at least two regular orbits on $V$.
\end{theorem}
\begin{proof}
This follows from \cite[Theorem 3.1]{YY2} and \cite[Theorem 3.1]{YY3}.
\end{proof}



\begin{lemma} \label{2.4}
Let $G$ be a finite solvable permutation group of degree $n$. Denote the number of cycles by $n(g)$ and the number of fixed points by $s(g)$. If $g \in G^{\#}$, then $n(g) \leq (n+s(g))/2 \leq (p+o(g)-1)n/(o(g) p) \leq 3n/4$.
\end{lemma}
\begin{proof}
Let $V$ be a minimal normal subgroup of $G$ and $S$ denote a point stabilizer. Then $n=|\Omega|=|V|$. If $s(g)=0$, we clearly have $n(g) \leq n/2$. We thus may assume that $g$ has fixed points, and without loss of generality $g \in S$. Since the actions of $S$ on $V$ and $\Omega$ are permutation isomorphic, it follows that $s(g)=|\bC_V(g)|$ and since $S$ acts faithfully on $V$, $s(g) \mid |V|/p =n/p$ where $|V|=p^n$. Therefore,

$n(g) \leq s(g)+(n-s(g))/o(g)=(n+(o(g)-1)s(g))/o(g) \leq (p+o(g)-1)n/(o(g) p) \leq \frac {3n} 4$.
\end{proof}

\begin{lemma}
Let $G$ be a finite solvable permutation group of degree $n$. Let $H$ be a nilpotent subgroup $G$, then we have $|H| \leq 2^n$.
\end{lemma}
\begin{proof}
\end{proof}



Since we need to compare the orbit size and the group order rather than to prove the existence of regular orbits, we need some quantitative results about how the solvable permutation group of odd order acts on the power set of the base set. The following result can be viewed as the quantitative analogue of a result of Gluck ~\cite[Corollary 5.7(b)]{manz/wolf} about primitive permutation groups.

\begin{theorem} \label{thmperm2}
Let $G$ be a solvable primitive permutation group of degree $n$ acts on a set $\Omega$ and $H$ be a nontrivial nilpotent subgroup of $G$, we show that if $k = \lceil \frac {\ln|H|}{\ln {2} } \rceil$, then for any subset $\Lambda$ of size $0 \leq m \leq k$, we may find  $\Delta \subseteq \Omega-\Lambda$ such that $|H:\stab_H(\Delta)|^2 \cdot 2^{m-1} \ge |H|$.
\end{theorem}
\begin{proof}
Assume $m=0$, then by $[7, \text{Corollary } 5.7]$, we can check that
\[|H:\stab_H(\Delta)|^2 \ge 2|H|\]
Assume $m > 1$.
Then, let $g$ be an element in $H$, we denote by $n(g)$ the number of cycles of $g$ on $\Omega$, and by $s(g)$ the number of fixed points.

For a subset $X \subseteq G$, it is worthwhile to consider in the following the set:

\[S(X)=\{(g,\Gamma) \ |\ g \in X, \Gamma \subseteq \Omega, g \in \stab_G(\Gamma) \}\]

Note that $g \in G$ stabilizes exactly $2^{n(g)}$ subsets of $\Omega$.
Since $H$ is a primitive solvable permutation group on the finite set $\Omega$, by $[7, \text{p. 39}]$, $H$ has a unique minimal normal subgroup $M$ with respect to a point stabilizer $H_\alpha$, where $\alpha \in \Omega$. From this, we also have that $H = H_\alpha \cdot M, H_\alpha \cap M = 1, C_H(M) = M$, and $M$ acts regularly on $\Omega$. Consequently, $|M| = |\Omega|$. So, we get that
\[|H| = |H_\alpha| \cdot |M| = |H_\alpha| \cdot |\Omega|\]
Combining this with by ~\cite[Theorem 3.3]{manz/wolf}, we get that
\[|H| \leq \frac{n^\beta}{2} \cdot n = \frac{n^{\beta+1}}{2}\]
Where $\beta = \frac{\log 32}{\log 9}$.

Therefore,
\[|S(H^{\#})| \leq \frac{n^{\beta+1}}{2}\cdot 2^{n(g)}\]




After taking away a subset $\Lambda$ from set $\Omega$, $|\Omega| - |\Lambda| = n - m$ elements remain Now, since $m \leq k$, the power set $P(\Omega - \Lambda)$ of this remaining set has size  $2^{n-m}$, and if $|S(H^\#)| < |P(\Omega - \Lambda)|$, then there must exist a regular orbit. So, finding a regular orbit is equivalent to finding when
\[|S(H^\#)| < 2^{n-m}\]
Which simplifies to
\[2^m \cdot \left(\frac{n^{\beta+1}}{2}\right) \cdot 2^{n(g)} < 2^n\]

If the nilpotent group $H$ has a regular orbit on $\Omega \oplus \Omega$, then Lemma 3.1 from \cite{manz/wolf} says that there is an orbit in $H$ where $|H:C_H(\Delta)|^2 \geq |H|$, $\Delta \in \Omega$.

Thus, if the following inequality is satisfied,
\[|H| \cdot (2^{n(g)})^2 \leq \left(\frac{2^n}{|H|}\right)^2\]
then $H$ has a regular orbit on $\Omega \oplus \Omega$.
Using Theorem 3.3 of \cite{manz/wolf}, we can bound $|H| \leq \frac{n^{\beta+1}}{2}$. Using Lemma \ref{2.4}, we can bound $n(g) \leq \frac{p+o(g)-1}{o(g)p}n \leq \frac{p+1}{2p}n \leq \frac{3n}{4}$ for $g \in G^\#$. Simplifying, we get
\[\frac{n^{\beta+1}}{2} \leq 2^\frac{n}{6}\]
Which holds for all $n \geq 97$

Here we split the remaining $n < 97$ into two cases: prime $n$ and prime power $n$.

First consider the prime $n=p$ case. If we now use the finer bound $n(g) \leq \frac{(p+1)n}{2p}$, we get the inequality
\[|H|^3 \cdot (2^\frac{(p+1)n}{2p})^2 \leq 2^{2n}\]
which simplifies to
\[n^3 \leq 2^{n-1}\]

which holds for all prime $n > 11$.

Now consider with the prime power cases (i.e. $n=p^\alpha, \alpha > 1$, we remain to check all such prime powers less than $97$. Using a GAP program, we obtain the following specific bounds for $|H|$:

\begin{table}[h!]
\centering
\begin{tabular}{|c|c|}
\hline
$n$ & largest $|H|$ \\
\hline
4  & 8   \\
8  & 8   \\
9  & 27  \\
16 & 128 \\
25 & 32  \\
27 & 81  \\
32 & 32  \\
49 & 96  \\
64 & 1024 \\
81 & 729 \\
\hline
\end{tabular}
\caption{}
\label{tab:nilpotent_selected}
\end{table}

Notice that, only $n=4, 8, 9, 16, 27$ fail to satisfy the inequality
\begin{equation}\label{bound}
|H|^3 \cdot (2^{n(g)})^2 \leq 2^{2n}
\end{equation}
When using the crude bound $n(g) \leq \frac{p+1}{2p}n$. So, to evaluate the smaller $n$, we run a GAP program to individually check that some subset $\Delta \subseteq \Omega - \Lambda $ satisfies the inequality in every nilpotent subgroup $H$ of all the solvable primitive permutation groups $G$ with degree $n$. This holds for primes $n = 2, 3, 5, 7, 11$ and prime powers $4, 8, 9, 16$ but takes too long to run for $n=27$. Instead, we will perform group analysis for this one remaining case.

A quick GAP program lists out all potential nilpotent $H \leq G$ for $\deg(G) = n = 27$. We find that the largest subgroup $H = (C_3 \times C_3 \times C_3) \rtimes C_3$ is the only subgroup of order $81$, and thus is the only subgroup that violates inequality \ref{bound}. All other $H$ satisfy $|H| \leq 27$ and satisfy \ref{bound}. Now, looking at $(C_3 \times C_3 \times C_3) \rtimes C_3$, we see that there is $1$ element of order $1$, $44$ elements of order $3$, and $36$ elements of order $9$. 

Now we use the finer bound $n(g) \leq \frac{p+o(g)-1}{o(g)p}n$, where $n = p^\alpha = 3^3 = 27$. For elements of order $o(g) = 3$, we have $n(g) \leq 15$. Clearly $n(g)$ for a higher order element will have a lower bound in this case. Therefore, showing that this same group with only order $3$ elements satisfies the inequality \ref{bound} implies that this inequality holds true for this group with some order $3$ and some order $9$ elements. Plugging into the inequality, we get

\[(81)^3 \cdot (2^{15})^2 \leq 2^{2\cdot27}\]

Which is true. Thus, we have verified the last remaining case $n=27$, and proving the theorem in its full generality. 
\end{proof}

\section{Main Theorems} \label{sec:maintheorem}

Isaacs proved the following result ~\cite[Theorem A]{IMI2}. Let $P$ be a nontrivial $p$-group that acts faithfully on a group $H$, where $|H|$ is not divisible by $p$. Then there exists an element $x \in H$ such that $|\bC_P(x)| \leq {(|P|/p)}^{1/p}$. By the Hartley-Turull's Lemma ~\cite[Lemma 2.6.2]{HTurull}, this is equivalent to the following statement: Let $P$ be a nontrivial $p$-group that acts faithfully on a vector space $V$, where $|V|$ is not divisible by $p$. Then there exists an element $x \in V$ such that $|\bC_P(x)| \leq (|P|/p)^{1/p}$. This compares well with the following Theorem.


\begin{theorem} \label{thm22}
  Suppose that $G$ is a finite solvable group and $V$ is a finite, faithful and completely reducible $G$-module. Assume $H$ is a nilpotent subgroup of $G$ and $2 \mid |H|$, then there exits $v \in V$ such that $2 \cdot |H| \leq |v^H|^2$.
\end{theorem}
\begin{proof}
We consider the minimal counterexample on $|G|+\dim V$.

First we show that we may assume that $V$ is irreducible. Assume not, then $V=V_1\oplus V_2$ where $V_1,V_2$ are nontrivial $G$-submodules of $V$.
Let $m_1$ be the largest orbit size of $H$ on $V_1$, and let $m_2$ be the largest orbit size in the action of $C_H(V_1)$ on $V_2$.
Moreover, let $v_i\in V_i$ ($i=1,2$) be representatives of these orbits. Put $v=v_1+v_2$. Then $C_H(v)=C_H(v_1)\cap C_H(v_2)$ and hence:
\begin{eqnarray*}
M&\geq&|v^H|=|H:(C_H(v_1)\cap C_H(v_2))|\\
 &=&|H:C_H(v_1)|\cdot|C_H(v_1):(C_H(v_1)\cap C_H(v_2))|\ =\ m_1\cdot|C_H(v_1):C_{C_H(v_1)}(v_2)|\\
 &=&m_1\cdot|v_2^{C_H(v_1)}|\  \geq\  m_1|v_2^{C_H(V_1)}|\ =\ m_1 m_2.
\end{eqnarray*}

Since $2 \mid |H|$, either $2 \mid |H/\bC_H(V_1)|$ or $2 \mid |\bC_H(V_1)|$. By the inductive hypothesis we further conclude that

\[(*)\quad M^{2}\geq m_1^{2} m_2^{2} \geq 2 \cdot |H/\bC_H(V_1)|\cdot |\bC_H(V_1)| = 2 \cdot |H|.\]

So from now on indeed let $V$ be irreducible.\\

We now assume that $V$ is not primitive. We hence assume that there exists a proper subgroup $L_1$ of $G$ and an irreducible $L_1$-submodule
$V_1$ of $V$ such that $V ={V_1}^G$. By transitivity of induction, we can choose $L_1$ to be a maximal subgroup of $G$.
In particular, $S \cong G/N$ is a primitive permutation group on a right transversal of $L_1$ in $G$, where $N$ is the normal
core of $L_1$ in $G$. Let $V_N=V_1 \oplus \cdots \oplus V_m$, where the $V_i$s are irreducible $L_i$-modules where $L_i=\bN_G(V_i)$ and
$m>1$. We know $G/N$ primitively permutes the elements of $\{V_1, \dots, V_m\}$. Clearly, $|H| = |HN/N| \cdot |H \cap N|$. \\

Define:
\[N_i=\bC_N\left(\sum_{j=1}^{i-1}V_j\right) \Bigg/  \bC_N\left(\sum_{j=1}^{i}V_j\right)\]
for $i=1,\dots ,m$ and note that $N_1=N/\bC_N(V_1)$ and $N_m=\bC_N(\sum_{j=1}^{m-1}V_j)$. We define $H_i$ to be the image of $H$ in $N_i$.
Then $|H \cap N| \leq \prod_{i=1}^{m} |H_i|$.
Clearly $N_i$ acts completely reducibly on $V_i$ for $i=1,\dots ,m$. Let $M_i$ be the largest orbit size of the action of $H_i$ on $V_i$
($i=1,\dots ,m$), and let $v_i\in V_i$ be representatives of the corresponding orbits for all $i$.
Thus, \[M_{H \cap N}\  \geq\ \left|\left(\sum_{i=1}^{m}v_i\right)^{H \cap N} \right| \ \geq\ \prod_{i=1}^{m} |v_i^{H_i}|\ =\ \prod_{i=1}^{m} M_i.\]

We define $\bar H=HN/N$ and the set of $m$ elements to be $\Omega$; we collect set the subset $\Omega$ that contains all $i$ such that $H_i>1$ to be $\Lambda$.

If $|\bar H| \leq 2^{|\Lambda|-1}$, then we set $w=w_1+w_2+\cdots+w_m$ where $w_i=v_i$ if $H_i>0$ and $w_i=0$ if $H_i=1$.\\

If $|\bar H| > 2^{|\Lambda|-1}$, then by Theorem ~\ref{thmperm2}, we may find $\Delta \subseteq \Omega-\Lambda$ such that $|\bar H:\stab_{\bar H}(\Delta)|^{2} \cdot 2^{|\Lambda|-1} \ge |\bar H|$. We set $w=w_1+w_2+\cdots+w_m$ where $w_i=v_i$ if $i \not\in \Delta$ and $w_i=0$ if $i \in \Delta$.\\

We further conclude that
\[\ M\ \geq\ |\bar H:\stab_{\bar H}(\Delta)| \cdot  \prod_{i \in \Lambda} M_i\]

We have that
\begin{equation}
\begin{aligned}
M^{2} \geq  |\bar H:\stab_{\bar H}(\Delta)|^{2} \cdot \prod_{i \in \Lambda} M_i^{2} \\ \geq  |\bar H:\stab_{\bar H}(\Delta)|^{2} \cdot 2^{|\Lambda|} \cdot \prod_{i \in \Lambda}|H_i|\ \\ \geq
2 \cdot |\bar H| \cdot |H \cap N| \\ =  2 |H|.
\end{aligned}
\end{equation}

Hence, now we may assume that $V$ is irreducible and primitive.

Step 2. If $V$ is primitive, then by the main results of ~\cite{YY2,YY3,YY7,Holt}, we know that either $e=1$ or $G$ has a regular orbit on $V$ except for a finite amount of cases listed in ~\cite{Holt}. For all those exceptional cases listed in ~\cite{Holt}, we may use GAP to verify that the result holds. 

Thus, we may assume $e=1$ and $G$ is a subgroup of $\Gamma(V)$. Since $H$ is a subgroups of $G$, we may set $|V|=p^n$ and $|H|=|H/C||C|$ where $|H/C| \mid n$ and $|C| \mid |V|-1$. For $g \in H \backslash C$, $|\bC_V(g)| \leq |V|^{1/2}$, and for $g \in C$, $|\bC_V(g)| \leq |V|^{1/2}$. If $|C| < |V|^{1/2}$, then $H$ has a regular orbit. Thus, we may assume that $|C|\geq |V|^{1/2}$. This shows that $|C|$ is pretty big. It suffices to show that $|C|^{2} \geq 2|H|$.

$|C| \geq 2 |H/C|$ and assume not, we have $|V|^{1/2} < 2 (\log_p |V|)$.

$|V|=2^2,2^3,3^2,2^4,3^3,2^5,2^6,2^7$.

 $|V|=2^2$. $|C|=3$, $|H/C|=2$, not nilpotent.

 $|V|=2^3$. $|C|=7$, $|H/C|=3$, not nilpotent.

 $|V|=3^2$. $|C|=8$, $|H/C|=2$, orbit of size 8.

 $|V|=2^4$. $|C|=15$, $|H/C|=4$, orbit of size $15$. $|C|=3$, $|H/C|=4$, regular orbit. $|C|=5$, $|H/C|=4$, not nilpotent.

 $|V|=3^3$. $|C|=13$, $|H/C|=3$, not nilpotent. $|C|=2$, $|H/C|=3$, regular orbit.

 $|V|=2^5$. $|C|=31$, $|H/C|=5$, orbit of size $31$.

 $|V|=2^6$. $|C|=9$, $|H/C|=6$, orbit of size $27$. $|C|=7$, $|H/C|=6$, orbit of size $21$.

 $|V|=2^7$. $|C|=127$, $|H/C|=7$, orbit of size $127$.
\end{proof}

\begin{theorem} \label{thm32}
Suppose that $G$ is a finite solvable group and $V$ is a finite, faithful and completely reducible $G$-module. Let $H$ be a nilpotent subgroup of $G$, then there exits $v \in V$ such that $|\bC_H(v)| \leq (|H|/p)^{1/p}$, where $p$ is the smallest prime divisor of $|H|$.
\end{theorem}
\begin{proof}

We first assume that $2 \mid |H|$, then the result follows by Theorem ~\ref{thm22}.

We now assume that $2 \nmid |H|$. Suppose $\CHA(\FF)=2$, then $H$ acts coprimely on $V$, and the result follows by the main result of ~\cite{IMI2}. Suppose that $\CHA(\FF) \neq 2$, then $H$ has a regular orbit on $V$ by the main result of ~\cite{YY5}.
\end{proof}

\begin{lemma}\label{lem4}
Let $G$ be a finite group, and suppose $N\unlhd G$ such that $G/N$ is nilpotent. Then there exists a nilpotent subgroup $U\leq G$ such that $G=NU$.
\end{lemma}
\begin{proof}
See, for example, \cite[III, Satz 3.10]{Huppert1}.
\end{proof}

The following result strengthens ~\cite[Theorem 1.2]{KellerYang}.

\begin{theorem}\label{thm33}
Let $G$ be a finite solvable group and $V$ a finite faithful completely reducible $G$-module, possibly of mixed characteristic.
Let $M$ be the largest orbit size in the action of $G$ on $V$. We denote $G^{\gN}$ to be the normal subgroup of $G$ such that $G/G^{\gN}$
is the maximum nilpotent quotient. We set $p$ to be the smallest prime dividing $|G/G^{\gN}|$. Then \[p \cdot |G/G^{\gN}|\leq M^{\frac p {p-1}}.\]
\end{theorem}
\begin{proof}
By Lemma \ref{lem4} there exists a nilpotent $L\leq G$ such that $G^{\gN}L=G$. Thus, there is a subgroup $H$ of $L$ such that $|G/G^{\gN}| \mid |H|$ and $\pi(|H|)=\pi(|G/G^{\gN}|)$ seems true, might need to replace $H$ by a smaller one (exclude the prime in $H \cap L$).

By Theorem \ref{thm32} it follows that if $M_H$ is the largest orbit size of $H$ on $V$, then $p|H|\leq M_H^{\frac p {p-1}}$. Hence, altogether $p|G/G^{\gN}|\leq p|H| \leq M_H^{\frac p {p-1}} \leq M^{\frac p {p-1}}$, as desired.
\end{proof}

How about central composition series. I think one can get the square bound anyway, this appeared in the paper with Qian.


\begin{prop}
Let $p$ be a given prime, let $\mathcal{F}$ be the class of abelian groups or the class of nilpotent groups. Assume that  $O_p(G) = 1$ and that  $G$ has an $\mathcal{F}$-quotient group of order $k$. Then $G$ has a  solvable subgroup $H$ with $O_p(H)=1$ such that $|H/H^{\mathcal{F}}|\geq k$.
\end{prop}

Can we still preserve the divisibility or not making $p$ smaller?

\begin{theorem}\label{thm34}
Let $G$ be a finite group and $V$ a finite faithful completely reducible $G$-module, possibly of mixed characteristic.
Let $M$ be the largest orbit size in the action of $G$ on $V$. We denote $G^{\gN}$ to be the normal subgroup of $G$ such that $G/G^{\gN}$
is the maximum nilpotent quotient. We set $p$ to be the smallest prime dividing $|G/G^{\gN}|$. Then \[|G/G^{\gN}|\leq M^{\frac p {p-1}}.\]
\end{theorem}

The following are some applications.

\begin{theorem} \label{thm2a}
Let $G$ be a finite group and $b(G)$ be the largest character degree of $G$, and we denote $H=G/F(G)$, then we have that $|H/H^{\gN}| \leq b^{\frac p {p-1}}(G)$ where $p$ is the smallest prime divisor of $|H/H^{\gN}|$.
\end{theorem}
\begin{proof}

Note that for any subgroup or quotient group $L$ of $G$, $b(L) \leq b(G)$. 


 Let $F^*(G)$ be the generalized Fitting subgroup of $G$ and we know that $F^*(G)=F(G)E(G)$. Let $Z=Z(E(G))$. By the proof of ~\cite[Lemma 3.15]{GLS2}
and by ~\cite[Lemma 3.16]{GLS2} we know that $F(G/Z)=F(G)/Z$, $E(G/Z)=E(G)/Z$, $F^*(G/Z)=F^*(G)/Z$,
and $C_{G/Z}(E(G/Z))=C_G(E(G))/Z$. Thus, we may assume $Z=1$.


 We will construct a character $\alpha \in \Irr(G)$ such that $\alpha^{\frac p {p-1}}(1) \geq |H/H^{\gN}|$.


 We first assume that $E(G)=1$. We know that $G/F(G)$ acts completely reducibly on $F(G)/\Phi(F(G))$.

 It follows by Theorem ~\ref{thm34} that there exists some linear $\lambda \in \Irr(F(G)/\Phi(F(G)))$ such that $b^{\frac p {p-1}}(G) \geq |H/H^{\gN}|$ and we choose $\alpha \in \Irr(G)$ such that $\alpha(1)=b(G)$. We have $\alpha^2(1) \geq |H/H^{\gN}|$. Therefore, the theorem holds in this case.

 We now assume that $E(G)>1$. Note that $E(G)$ is a direct product of simple groups, and let $E(G)=E_1 \times E_2 \cdots \times E_m$ where $E_i$ is a direct product of $k_i$ isomorphic simple groups $L_i$.

 We first consider $\bC_G(E_m)$, and we denote $N_m=\bC_G(E_m)$ and $G_m=G/N_m$. Since $E_m \cong (E_m \times N_m)/N_m$, we can consider $E_m$ as a subgroup of $G_m$, and we know that $G_m$ is embedded in $\Aut(E_m) \cong \Aut(L_m) \wr S_{k_m}$.

 Recursively, we define $N_{i}=\bC_{N_{i+1}}(E_i)$, $G_i=N_{i+1}/N_i$ for $1 \leq i \leq m-1$. Since $E_i \cong (E_i \times N_i)/N_i$, we can consider $E_i$ as a subgroup of $G_i$, and we know that $G_i$ is embedded in $\Aut(E_i) \cong \Aut(L_i) \wr S_{k_i}$.

 So we have a chain of quotients, and we know that $N_1=\bC_{N_2}(E_{1})=\bC_G(E(G))$. Thus, $|H/H^{\gN}| \mid \prod_{i=1}^{m} |G_i/G_i^{\gN}| \cdot |N_1/N_1^{\gN}|.$

 Clearly $F(G) \leq N_1$ and $N_1 \cap E(G)=1$. By ~\cite[Lemma 3.13(i)]{GLS2}, $F^*(N_1)=F(N_1)=F(G)$. We have $\bC_{N_1}(F(G)) \leq F(G)$. We know that $N_1/F(G)$ acts faithfully and completely reducibly on $F(G)/\Phi(F(G))$.


 If $H_0=1$, then clearly $b(N_1) \geq \sqrt{|H_0|}$. If $H_0 \neq 1$, it follows by Theorem ~\ref{t206} that there exists some linear $\lambda \in \Irr(F(G)/\Phi(F(G)))$ such that $b(N_1) \geq \sqrt{|N_1/N_1^{\gN}|}$ and we choose $\psi \in \Irr(N_1)$ such that $\psi(1)=b(N_1)$. We have $\psi^2(1) \geq |N_1/N_1^{\gN}|$ and the inequality is strict if $H_0>1$. 

Let $\chi=\chi_1 \times \cdots \times \chi_m$ and we have that $\chi^2(1) \geq \prod_{i=1}^{m} |G_i/G_i^{\gN}|$.

Note that $\alpha=\chi \times \psi$ is an irreducible character of $E(G) \times N_1$. It follows that $b^2(G) \geq  b^2(E(G) \times N_1) \geq \alpha^2(1) = \chi^2(1) \psi^2(1) \geq |H|$.
\end{proof}

We strengthen a result of Qian and Yang ~\cite[Theorem 1.3]{YY16} (the nilpotent case).

\begin{theorem} \label{thm2b}
Let a finite group $G$ act faithfully on a finite group $V$,
and $M$ be the largest orbit size in the action of $G$ on
$V$. Then any one of the following conditions guarantees that
$$ M^2 \geq 2 \cdot 2^{\mu(G)}\cdot {\rm Ord}_{\mathcal{N}}(G).$$

{\rm (1)} $V$ is a $p$-group and $O_p(G)=1$ for some prime $p$;

{\rm (2)} $V$ is  a completely reducible $G$-module, possibly of mixed characteristic.
\end{theorem}

The following is another result about nilpotent subgroups of solvable groups.

\begin{theorem} \label{thm2c}
Let $G$ be a finite solvable group and $b(G)$ be the largest character degree of $G$. Let $H$ to be a nilpotent subgroup of $G$, then we have that $p|H F(G)/F(G)| \leq b^{\frac p {p-1}}(G)$ where $p$ is the smallest prime divisor of $|H F(G)/F(G)|$.
\end{theorem}

\begin{theorem} \label{app1}
  Suppose that $G$ is a finite solvable group and let $H$ be a nilpotent Hall $\pi$-subgroup, then $p|G:O_{\pi' \pi}(G)|_{\pi} \leq b(H)^{\frac p {p-1}}$ where $p$ is the smallest number in $\pi$.
\end{theorem}
\begin{proof}
We may assume that $O_{\pi'}(G)=1$. Let $N=O_{\pi}(G)$. Then, fairly standard arguments show that $C=\bC_G(\bF(N)/\Phi(N))\subseteq N$. Write $V=\Irr(\bF(N)/\Phi(N))$ and $\bar{G}=G/C$. Thus $O_{\pi}(\bar{G})=N/C$. Now, $V$ is a faithful $\bar{G}$-module such that $V_{O_{\pi}(\bar{G})}$ is completely reducible.

Let $H$ be a nilpotent Hall $\pi$-subgroup of $G$ and let $\bar{H}=H/N$. By Theorem ~\ref{thm32}, there exists $\lambda \in V$ such that $|\bC_{\bar{H}}(\lambda)| \leq (|\bar{H}|/p)^{1/p}$. Let $\xi \in \Irr(\bC_{H}(\lambda)|\lambda)$ and $\alpha=\xi^H \in \Irr(H)$. Thus, $p|G:N|_{\pi} \leq p|\bar{H}| \leq \alpha(1)^{\frac p {p-1}} \leq b(H)^{\frac p {p-1}}$, as wanted.
\end{proof}
Remark: This strengthens ~\cite[Theorem 4.2]{YYsubgroup}.



\begin{cor} \label{app3}
  Suppose that $G$ is a finite solvable group and let $P \in \Syl_2(G)$, then $2|G:O_{2'2}(G)|_2 \leq b(P)^2$.
\end{cor}
\begin{proof}
This follows from Theorem ~\ref{app1} by choosing $H$ to be a Sylow $2$-subgroup of $G$.
\end{proof}

\section{Acknowledgement} \label{sec:Acknowledgement}

This work was partially supported a grant from the Simons Foundation (\#918096, to YY).







\end{document}